\documentstyle[11pt,amstex,amssymb]{amsart}\textheight=19.8truecm\textwidth=14.1truecm
\begin{document}

\centerline{\Large \bf Matrix Subadditivity Inequalities and
Block-Matrices}

\vskip 25pt \centerline{{\large \sl \ Jean-Christophe Bourin }}

\vskip 35pt \noindent {\small {\bf Abstract.} We give a number of
subadditivity results and conjectures for symmetric norms, matrices
and block-matrices. Let $A$, $B$, $Z$ be
 matrices of same size and suppose that $A,\, B$ are normal and $Z$ is expansive, i.e., $Z^*Z\ge
I$. We conjecture that
$$
\Vert\, f(|A+B|)\,\Vert \le \Vert\, f(|A|)+f(|B|)\,\Vert \quad {\rm
and} \quad \Vert\, f(|Z^*AZ|)\,\Vert \le \Vert\, Z^*f(|A|)Z\,\Vert
$$
for all non-negative concave function $f$ on $[0,\infty)$  and all
symmetric norms $\|\cdot\|$ (in particular for all Schatten
$p$-norms). This would extend known results for positive operator to
all normal operators. We prove these inequalities in several cases
and we propose some related open questions, both in the positive and
normal cases. As nice applications of subadditivity results we get
some unusual estimates for partitioned matrices. For instance, for
all symmetric norms and $0\le p\le 1$,
$$
\left\| \,\left|
\begin{pmatrix}A&B \\C&D \end{pmatrix}
 \right|^p \,\right\|
  \le \left\|\, |A|^p+|B|^p+|C|^p+|D|^p\,\right\|.
 $$
whenever the partitioned matrix is Hermitian or its entries are
normal. We conjecture that this estimate for $f(t)=t^p$ remains true
for all non-negative concave functions $f$ on the positive
half-line. Some results for general block-matrices are also given.


Mathematical subjects classification:   15A60, 47A30, 47A60}

\vskip 25pt\noindent {\large\bf  Introduction} \vskip 10pt

This paper complements  three recent works  on subadditivity type
inequalities for concave functions of positive operators
 \cite{B2},
\cite{BL} and \cite{BU} . These inequalities are matrix versions of
obvious scalars
 inequalities like
 $$
 f(za)\le zf(a) \quad {\rm and} \quad f(a+b) \le f(a)+f(b)
 $$
 for  non-negative concave functions $f$ on $[0,\infty)$  and
 scalars $a,\,b\ge0$ and $z\ge1$. By matrix version we mean
 suitable extensions where scalars $a,\,b,\dots z$ are replaced by $n$-by-$n$ complex
 matrices, i.e., operators on an
  $n$-dimensional Hilbert space ${\mathcal{H}}$, $A,\,B,\dots Z$.  For instance, if $A$ is positive (semi-definite) and $Z$ is
  expansive we know from \cite{B1}  a remarkable trace inequality, companion to
  one of the above scalar inequality,
    \begin{equation}
   {\mathrm Tr\,}f(Z^*AZ) \le {\mathrm Tr\,}Z^*f(A)Z.
   \end{equation}
This trace inequality may   be generalized
  by using the class of
 symmetric (or unitarily invariant)
 norms. Such
 norms satisfy  $\Vert A\Vert = \Vert UAV\Vert$ for all $A$ and all
unitaries $U,\, V$. In the next section we aim to extend  known
inequalities for the cone of positive operators to the set of all
normal operators. We have only partial results and  some open
problems are considered. A basic one would be to know if (1) can be
extended to normal operators $A$ as follows:
\begin{equation}
   {\mathrm Tr\,}f(|Z^*AZ|) \le {\mathrm Tr\,}Z^*f(|A|)Z.
   \end{equation}
The main part of this paper is Section 2 where we show how some
subadditivity results entail several new estimates for block
matrices. A special case involving four normal operators $A,\, B,\,
C,\, D$ is
\begin{equation*}
\left\| \,\left|
\begin{pmatrix}A&D \\C&B \end{pmatrix}
 \right|^p \,\right\|
  \le \left\|\, |A|^p+|B|^p+|C|^p+|D|^p\,\right\|
\end{equation*}
for all symmetric norms and $0\le p\le 1$. Such an inequality also
 holds in the important case of a Hermitian partitioned matrix. These estimates differ
from the usual ones in the literature where the norm of the full
matrix is evaluated with  the norms of its blocks, for instance see
\cite{Aud}, \cite{Kg}, \cite{KgNa} and \cite{BhKi1}. Our results
start the study of the general problem of comparing an operator on
${\mathcal{H}}\oplus {\mathcal{H}}$ with an operator on
${\mathcal{H}}$, more precisely, comparing the block-matrix
expression
\begin{equation*}
f\left(\left|
\begin{pmatrix}A&D \\C&B \end{pmatrix}
 \right|\right)
\end{equation*}
with the sum
$$
f(|A|)+f(|B|)+f(|C|)+f(|D|).
$$
Some natural conjectures, parallel to those ones of Section 1, are
naturally proposed.

The proofs use standard two by two block-matrix technics and basic
facts on majorisation, log-majorisation, convex functions and norms.
This background can be found in text books like \cite{Bh} and will
be used  without further reference.

 \vskip 20pt\noindent {\large\bf 1. Some recent and new subadditivity results}

\vskip 10pt
 If $A$ is positive (semi-definite), resp.\ positive
definite,
  we write $A\ge 0$, resp.\ $A>0$. If  $Z^*Z$ dominates the identity
  $I$, we say that $Z$ is expansive. The trace inequality (1) is a
  special case of the following theorem \cite{B2}:

 \vskip 10pt

\vskip 10pt\noindent
 {\bf Theorem 1.1a.} {\it Let $A\ge0$ and let $Z$  be expansive. If
 $f:[0,\infty)\to [0,\infty)$ is concave, then, for all symmetric norms,
 $$\Vert\, f(Z^*AZ)\,\Vert \le \Vert\, Z^*f(A)Z\,\Vert.$$
}

\vskip 10pt\noindent This estimate for congruences have been
completed in \cite{BU} with an estimate for sums:
 \vskip 10pt\noindent
 {\bf Theorem 1.1b.} {\it Let $A,\,B\ge0$ and let
 $f:[0,\infty)\to [0,\infty)$ be concave. Then, for all symmetric norms,
 $$\Vert\, f(A+B)\,\Vert \le \Vert\, f(A)+f(B)\,\Vert.$$
 }

\vskip 10pt\noindent This result  closed a long list of papers by
several authors, including Ando-Zhan \cite{AZ} and Kosem \cite{Ko},
on norm subadditivity inequalities,  cf.\ \cite{BL}, \cite{BU} and
references therein.
 Theorems 1.1a and 1.1b are nicely compatible. There are two aspects of a unique statement \cite{BL}:

 \vskip 10pt\noindent
 {\bf Theorem 1.1.}  {\it Let $\{A_i\}_{i=1}^m$ be positive, let  $\{Z_i\}_{i=1}^m$ be expansive and let
 $f$ be a non-negative concave function on $[0,\infty)$. Then, for all symmetric norms,
$$
\left\|f\left(\sum Z_i^*A_iZ_i\right) \right\| \le \left\| \sum
Z_i^*f(A_i)Z_i \right\|.
$$
}

\vskip 10pt\noindent It would be elegant and interesting to state
these theorems in the more general  framework of positive linear
maps $\Phi$ between matrix algebras. Indeed the conclusions of these
theorems look like statements for positive linear maps. For
instance, with
$$
\Phi\left(
\begin{pmatrix}
S &X\\
Y&T \end{pmatrix} \right) =S+T \qquad {\mathrm{and}}\qquad
M=\begin{pmatrix}
A &0\\
0&B \end{pmatrix} $$
 Theorem 1.1b claims
 \begin{equation*}
 \Vert f\circ\Phi(M)\Vert \le \Vert \Phi\circ f(M)\Vert.
 \end{equation*}
 Therefore we ask two questions.

\vskip 10pt\noindent
 {\bf Question 1.} Let $\Phi:{\mathrm{M}}_n\to{\mathrm{M}}_k$ be a
 positive linear map between the algebras of matrices of size $n$ and $k$ and let
 ${\mathcal{F}}_n$ and ${\mathcal{F}}_k$ be the sets of non-expansive, positive operators
 in ${\mathrm{M}}_n$ and ${\mathrm{M}}_k$ respectively. Does
 $$
 \Phi^{-1}({\mathcal{F}}_k) \subset {\mathcal{F}}_n
  $$
imply
\begin{equation*}
 \Vert f\circ\Phi(M)\Vert \le \Vert \Phi\circ f(M)\Vert
 \end{equation*}
 for all $M\ge0$, all symmetric norms and all non-negative concave
 functions on the positive half line ?

 \vskip 10pt\noindent
 {\bf Question 2.} Let $l: \Phi({\mathrm{M}}_n)\to{\mathrm{M}}_n$ be a lifting
 for $\Phi$, i.e.,\ an injective map such that $\Phi\circ l$ is the
 identity mapping of $\Phi({\mathrm{M}}_n)$. Let
  ${\mathcal{E}}_k$ be the set of expansive, positive operators
 in  ${\mathrm{M}}_k$. Does
 $$
 \Phi\circ l({\mathcal{E}}_k) \subset {\mathcal{E}}_k
  $$
imply
\begin{equation*}
 \Vert f\circ\Phi(M)\Vert \le \Vert \Phi\circ f(M)\Vert
 \end{equation*}
 for all $M\ge0$, all symmetric norms and all non-negative concave
 functions on the positive half line ?

\vskip 10pt
 Now we turn to some consequences of Theorem 1.1. The usual  operator norm
 is denoted by $\Vert \cdot \Vert_{\infty}$.

\vskip 10pt\noindent
 {\bf Corollary 1.2.}  {\it Let $\{A_i\}_{i=1}^m$ be normal, let  $\{Z_i\}_{i=1}^m$ be expansive and let
 $f$ be a non-negative concave function on $[0,\infty)$. Then,
$$
\left\|f\left(\left|\sum Z_i^*A_iZ_i\right|\right) \right\|_{\infty}
\le \left\| \sum Z_i^*f(|A_i|)Z_i \right\|_{\infty}.
$$
}

\vskip 10pt\noindent
 {\bf Proof.} Note that
 $$
 \begin{pmatrix}
|A_i| &A_i \\
 A^*_i &|A_i| \ \end{pmatrix}\ge 0
  $$
hence $$
\begin{pmatrix}
\sum Z^*_i|A_i|Z_i & \sum Z^*_iA_iZ_i \\
 \sum Z^*_iA^*_iZ_i & \sum Z^*_i|A_i|Z_i \ \end{pmatrix}\ge 0
 $$
so that
 $$\left\| \sum Z^*_iA_iZ_i\right\|_{\infty} \le \left\| \sum
 Z^*_i|A_i|Z_i\right\|_{\infty}.$$
 Since $f$ is non-decreasing and non-negative we have
 \begin{align*} \left\|f\left(\left|\sum Z_i^*A_iZ_i\right|\right)
 \right\|_{\infty} &=
 f\left(\left\|\sum Z_i^*A_iZ_i\right\|_{\infty}\right) \\
 &\le f\left(\left\|\sum Z_i^*|A_i|Z_i\right\|_{\infty}\right) \\
&=\left\| f\left(\sum Z_i^*|A_i|Z_i \right)\right\|_{\infty}.
\end{align*}
We then apply Theorem 1.1 to get the conclusion. \qquad $\Box$

\vskip 10pt\noindent We may propose:
 \vskip 10pt\noindent
 {\bf Conjecture 1.}  Corollary 1.2 holds for all symmetric norms.
 In particular the trace inequality (2) holds.

\vskip 10pt\noindent A special case of Corollary 1.2 is a matrix
version of an obvious inequality for complex numbers $z=a+ib$.\vskip
10pt\noindent
 {\bf Corollary 1.3.}  {\it Let   $Z=A+iB$ be a  decomposition in real and imaginary parts, and let
 $f$ be a non-negative concave function on $[0,\infty)$. Then,
$$
\left\|f(|Z|) \right\|_{\infty} \le
\left\|f(|A|)+f(|B|)\right\|_{\infty}.
$$
}

\vskip 10pt\noindent The next application of Theorem 1.1 gives a
partial answer to Conjecture 1. We first give a simple definition.

\vskip 10pt\noindent {\bf Definition.} A function $f(t)$  on
$[0,\infty)$ is {\bf e-convex} if $f(e^t)$ is convex on
$(-\infty,\infty)$. Some non-negative concave but e-convex functions
are $f(t)=t^p$, $0\le p\le 1$, and $f(t) =\log(1+t)$.

\vskip 15pt\noindent
 {\bf Corollary 1.4.}  {\it Let $\{A_i\}_{i=1}^m$ be normal, let  $\{Z_i\}_{i=1}^m$ be expansive and let
 $f$ be a non-negative concave and e-convex function on $[0,\infty)$.
 Then, for all symmetric norms,
$$
\left\|f\left(\left|\sum Z_i^*A_iZ_i\right|\right) \right\| \le
\left\| \sum Z_i^*f(|A_i|)Z_i \right\|.
$$
}

 \vskip 10pt\noindent
Thus we have the following two special cases:
 \vskip 15pt\noindent
 {\bf Corollary 1.5.}  {\it Let $A$, $B$ be normal and let
 $f$ be a non-negative concave and e-convex function on $[0,\infty)$.
 Then, for all symmetric norms,
$$
\Vert f(|A+B|)\Vert \le \Vert \,f(|A|)+f(|B|)\,\Vert
$$
}

\vskip 10pt\noindent
 {\bf Corollary 1.6.}  {\it Let $A$, $B$ be normal. Then, for all symmetric norms and $0\le p\le1$,
$$
\Vert \,|A+B|^p\,\Vert \le \Vert\, |A|^p+|B|^p\,\Vert
$$
}

\vskip 10pt\noindent We turn to the proof of Corollary 1.4.  Given
$A,\, B\ge 0$ we use the following notations. The weak-log
majorisation relation
$$
A\prec_{wlog} B
$$
means that for all $k=1,\, 2,\ldots$ we have
$$
\prod_{j=1}^k \lambda_j(A)\le \prod_{j=1}^k \lambda_j(B)
$$
where $\lambda_j(\cdot)$ are the eigenvalues arranged in decreasing
order and repeated according their multiplicities. The weak
majorisation relation
$$
A\prec_{w} B
$$
means that for all $k=1,\, 2,\ldots$ we have
$$
\sum_{j=1}^k \lambda_j(A)\le \sum_{j=1}^k \lambda_j(B).
$$

\vskip 10pt\noindent
 {\bf Proof of Corollary 1.4.} As in the proof of Corollary 1.2, we have
 $$
\begin{pmatrix}
\sum Z^*_i|A_i|Z_i & \sum Z^*_iA_iZ_i \\
 \sum Z^*_iA_iZ_i & \sum Z^*_i|A_i|Z_i \ \end{pmatrix}\ge 0
 $$
 so that for some contraction $K$
 $$
\sum Z^*_iA_iZ_i=\left(\sum Z^*_i|A_i|Z_i \right)^{1/2}K\left(\sum
Z^*_i|A_i|Z_i \right)^{1/2}
 $$
which implies via Horn's inequalities
$$
\left|\sum Z^*_iA_iZ_i\right|  \prec_{wlog}\sum Z^*_i|A_i|Z_i.
$$
Therefore, since a non-decreasing convex function preserves
weak-majorisation, the e-convexity property of $f$ entails
$$
f\left(\left|\sum Z^*_iA_iZ_i\right|\right)  \prec_{w}f\left(\sum
Z^*_i|A_i|Z_i\right) .
$$
Since $f$ is non-negative, this statement is equivalent to
$$
\left\|f\left(\left|\sum Z^*_iA_iZ_i\right|\right)\right\| \le
\left\|f\left(\sum Z^*_i|A_i|Z_i\right)\right\|
$$
for all symmetric norms. Hence,  Theorem 1.1 shows the result.
\qquad $\Box$

\vskip 15pt\noindent Theorem 1.1 also yields some results for
non-negative convex functions vanishing at $0$, see \cite{BL}. For
instance:
 \vskip 15pt\noindent
 {\bf Corollary 1.7.}  {\it Let $\{A_i\}_{i=1}^m$ be positive, let  $\{Z_i\}_{i=1}^m$ be expansive and  let $p\ge1$ .
  Then, for all symmetric norms,
$$
\left\| \sum Z_i^*A_i^pZ_i \right\| \le \left\|  \left(\sum
Z_i^*A_iZ_i \right)^p \right\|.
$$
}

\vskip 10pt\noindent The important special case
$$
||A^{p}+B^{p}||\leq ||(A+B)^p||
$$
for positive $A,\,B$ being due to Ando-Zhan \cite{AZ} and  to
Bhatia-Kittaneh \cite{BhKi2} for integer exponents. This inequality
raises some questions.

\vskip 10pt\noindent {\bf Question 3.} Given  $A,\,B\ge 0$ and
$p,\,q\ge 0$, is it true that
$$||A^{p+q}+B^{p+q}||\leq
||(A^p+B^p)(A^q+B^q)||\ ?$$ More generally, let $p=\sum p_i$ and
$q=\sum q_i$ where $p_i,\,q_i\ge 0$, $i=1,\cdots,m$. Does
$$
||A^p+B^p||\leq ||(A^{p_1}+B^{q_1})\cdots (A^{p_n}+B^{q_n})||
$$
hold ? One may also ask wether stronger inequalites like
$$||A^{p+q}+B^{p+q}||\leq
||(A^p+B^p)^{1/2}(A^q+B^q)(A^p+B^p)^{1/2}||$$ hold true.

 \vskip 10pt\noindent {\bf Question 4.} Given  $A,\,B\ge 0$ and
$p,\,q\ge 0$, is it true that
$$||A^pB^q+B^pA^q||\leq ||A^{p+q}+B^{p+q}||\ ?$$
Note that the similar inequality for the Heinz means
$$||A^pB^q+A^qB^p||\leq ||A^{p+q}+B^{p+q}||$$
is well-known. See \cite{HKos} for much more on matrix means.

\vskip 20pt\noindent {\large\bf 2. Applications to block-matrices }

\vskip 10pt Let ${\mathbb{A}}=[A_{i,\,j}]$ be an arbitrary
block-matrix and let $\Vert\cdot\Vert_1$ stand for the trace norm.
Then, for all non-negative concave function $f$ on $[0,\infty)$,
\begin{equation}
\Vert\, f(|{\mathbb{A}}|)\,\Vert_1 \le  \sum \left\|\,
f(|A_{i,\,j}|) \,\right\|_1.
\end{equation}
Indeed, this follows from Rotfeld's inequality \cite{Ro}: {\it For
arbitrary operators $X,\,Y$, we have ${\mathrm{Tr}}\, f(|X+Y|)\le
{\mathrm{Tr}}\, f(|X|) +{\mathrm{Tr}}\, f(|Y|) $}. As a consequence
of Theorem 1.1b we state a result for all symmetric norms:

 \vskip 10pt\noindent
 {\bf Theorem 2.1.}  {\it Let ${\mathbb{A}}=[A_{i,\,j}]$ be a block matrix with normal
 entries and let
 $f$ be a non-negative concave,  e-convex function on $[0,\infty)$.
  Then, for all symmetric norms,
$$
\left\|\, f(|{\mathbb{A}}|) \,\right\| \le \left\|\,\sum
f(|A_{i,\,j}|) \,\right\|.
$$
}

\vskip 10pt\noindent For $f(t)=t$ we recapture (with a simpler
proof) an observation from \cite{BU},
$$
\left\| {\mathbb{A}} \right\| \le \left\|\,\sum |A_{i,\,j}|
\,\right\|.
$$
Applying a non-negative concave function $f(t)$ to both sides in the
operator norm case and using Theorem 1.1b we then obtain:

 \vskip 10pt\noindent
 {\bf Corollary 2.2.}  {\it Let ${\mathbb{A}}=[A_{i,\,j}]$ be a block matrix with
 normal entries and let
 $f$ be a non-negative concave function on $[0,\infty)$. Then,

$$
\left\|\, f\left(\left|{\mathbb{A}}\right|\right)
\,\right\|_{\infty} \le \left\|\,\sum f(|A_{i,\,j}|)
\,\right\|_{\infty}.
$$
}

\vskip 10pt\noindent This result combined with (3) strongly
suggests:

 \vskip 10pt\noindent
 {\bf Conjecture 2.}  Corollary 2.2 holds for all symmetric norms.

 \vskip 10pt\noindent
Theorem 2.1 partially answers this conjecture. We turn to the proof
of Theorem 2.1. We start with two elementary, well-known lemmas.

\vskip 10pt\noindent
 {\bf Lemma 1.} {\it Let $A,\,B,\,X,\,Y \ge 0$ such that $B\prec_w Y$
 and $A\prec_w X$. Then,
$$
\begin{pmatrix}A&0 \\0&B \end{pmatrix}\prec_w \begin{pmatrix}X&0 \\0&Y
\end{pmatrix}.
$$
}

\vskip 10pt\noindent
 {\bf Proof.} We have
 $$
 \sum_{j=1}^k\lambda_j(A\oplus B)=\max_{s+t=k}\left\{\sum_{j=1}^s\lambda_j(A) +
 \sum_{j=1}^t\lambda_j(B)\right\}.
 $$
  Combining this with
 $$
 \sum_{j=1}^s\lambda_j(A) + \sum_{j=1}^t\lambda_j(B)\le  \sum_{j=1}^s\lambda_j(X) + \sum_{j=1}^t\lambda_j(Y)\le  \sum_{j=1}^k\lambda_j(X\oplus Y)
 $$
 ends the proof. \qquad $\Box$

\vskip 10pt\noindent
 {\bf Lemma 2.} {\it Let $A,\,B \ge 0$. Then,
$$
\begin{pmatrix}A&0 \\0&B \end{pmatrix}\prec_w \begin{pmatrix}A+B&0
\\0&0
\end{pmatrix}.
$$
}

\vskip 10pt\noindent {\bf Proof.} Note that
$$
\begin{pmatrix}A+B&0\\0&0
\end{pmatrix}=
\begin{pmatrix}A^{1/2}&B^{1/2}\\0&0
\end{pmatrix}
\begin{pmatrix}A^{1/2}&0\\B^{1/2}&0
\end{pmatrix}
$$
so that
$$
\begin{pmatrix}A+B&0\\0&0
\end{pmatrix}
\simeq
\begin{pmatrix}A&A^{1/2}B^{1/2}\\B^{1/2}A^{1/2}&B
\end{pmatrix}
\simeq
\begin{pmatrix}A&-A^{1/2}B^{1/2}\\-B^{1/2}A^{1/2}&B
\end{pmatrix}
$$
where $\simeq$ means unitarily congruent. Combining with
$$
\begin{pmatrix}A&0 \\0&B \end{pmatrix}
=\frac{1}{2}\begin{pmatrix}A&A^{1/2}B^{1/2}\\B^{1/2}A^{1/2}&B\end{pmatrix}
+\frac{1}{2}\begin{pmatrix}A&-A^{1/2}B^{1/2}\\-B^{1/2}A^{1/2}&B\end{pmatrix}
$$
gives the lemma. \qquad $\Box$

 \vskip 10pt\noindent {\bf Proof of
Theorem 2.1.}  We prove the theorem   for a partition in four normal
blocks
$$
{\mathbb{A}}=\begin{pmatrix}S&R \\T&Q \end{pmatrix}
$$
the proof of a more general partition being quite similar. Let
$$
\tilde{{\mathbb{A}}}=\begin{pmatrix}0&{\mathbb{A}}
\\{\mathbb{A}}^*&0
\end{pmatrix}
$$
and note that
\begin{equation*}
|\tilde{{\mathbb{A}}}|=\begin{pmatrix}|{\mathbb{A}}^*|&0
\\ 0&|{\mathbb{A}}|\end{pmatrix}
\end{equation*}
so that
\begin{equation}
|\tilde{{\mathbb{A}}}|\simeq\begin{pmatrix}|{\mathbb{A}}|&0
\\ 0&|{\mathbb{A}}|\end{pmatrix}
\end{equation}
where the symbol $\simeq$ stands for unitarily equivalent. On the
other hand
 $$
 \tilde{{\mathbb{A}}}=\tilde{{\mathbb{S}}}+\tilde{{\mathbb{T}}}+\tilde{{\mathbb{R}}}+\tilde{{\mathbb{Q}}}
 $$
 where
 $$
\tilde{{\mathbb{S}}}=\begin{pmatrix} 0&0&S&0 \\ 0&0&0&0
\\ S^*&0&0&0 \\ 0&0&0&0\end{pmatrix}
\quad \tilde{{\mathbb{T}}}=\begin{pmatrix} 0&0&0&0 \\ 0&0&T&0
\\ 0&T^*&0&0 \\ 0&0&0&0\end{pmatrix}
$$
and $$ \quad \tilde{{\mathbb{R}}}=\begin{pmatrix} 0&0&0&R \\ 0&0&0&0
\\ 0&0&0&0 \\ R^*&0&0&0\end{pmatrix}
\quad \tilde{{\mathbb{Q}}}=\begin{pmatrix} 0&0&0&0 \\ 0&0&0&Q
\\ 0&0&0&0 \\ 0&Q^*&0&0\end{pmatrix}
 $$
 are Hermitian. Therefore, as in the proof of Corollary 1.4,
 $$
 |\tilde{{\mathbb{A}}}|\prec_{wlog}
 |\tilde{{\mathbb{S}}}|+|\tilde{{\mathbb{T}}}|+|\tilde{{\mathbb{R}}}|+|\tilde{{\mathbb{Q}}}|,
 $$
 that is
 $$
 |\tilde{{\mathbb{A}}}|\prec_{wlog}
 \begin{pmatrix} |S^*|+|R^*|&0&0&0 \\ 0&|T^*|+|Q^*|&0&0
\\ 0&0&|S|+|T|&0 \\ 0&0&0&|R|+|Q|\end{pmatrix}.
$$
Since $f$ is non-decreasing and e-convex we infer
 $$
 f(|\tilde{{\mathbb{A}}}|)\prec_{w}
 \begin{pmatrix} f(|S^*|+|R^*|)&0&0&0 \\ 0&f(|T^*|+|Q^*|)&0&0
\\ 0&0&f(|S|+|T|)&0 \\ 0&0&0&f(|R|+|Q|)\end{pmatrix}.
$$
By Theorem 1.1b and Lemma 1 we then get
 $$
 f(|\tilde{{\mathbb{A}}}|)\prec_{w}
 \begin{pmatrix} f(|S^*|)+f(|R^*|)&0&0&0 \\ 0&f(|T^*|)+f(|Q^*|)&0&0
\\ 0&0&f(|S|)+f(|T|)&0 \\ 0&0&0&f(|R|)+f(|Q|)\end{pmatrix}.
$$
Gathering the two first lines, and the two last ones, we have via
Lemmas 2 and 1
$$
 f(|\tilde{{\mathbb{A}}}|)\prec_{w}
 \begin{pmatrix} f(|S^*|)+f(|T^*|)+f(|R^*|)+f(|Q^*|)&0\\ 0&f(|S|)+f(|T|)+f(|R|)+f(|Q|)\end{pmatrix}.
$$
By using (4) we then obtain, using normality of $S,\,T,\,R,\,Q$,
$$
f(|{\mathbb{A}}|)\prec_{w} f(|S|)+f(|T|)+f(|R|)+f(|Q|)
$$
which is equivalent to inequalities for symmetric norms. \qquad
$\Box$

\vskip 10pt Let us point out two variations of Theorem 2.1 in which
some operators are not necessarily normal.

\vskip 10pt\noindent
 {\bf Corollary 2.3.} {\it Let ${\mathbb{T}}$ be a triangular block-matrix
 $$
 {\mathbb{T}}=\begin{pmatrix} A&N\\0&B\end{pmatrix}.
 $$
 in which $N$ is normal. Let $f:[0,\infty)\to[0,\infty)$ be  concave and e-convex. Then, for all symmetric
 norms,
 $$
\Vert f(|{\mathbb{T}}|)\Vert \le \Vert f(|A^*|)+ f(|N|)+f(|B|)
\Vert.
 $$
 }

\vskip 10pt\noindent
 {\bf Proof.} Consider the polar decompositions $A=|A^*|U$ and
 $B=V|B|$, note that
 $$
\left|\begin{pmatrix} A&N\\0&B\end{pmatrix}\right|\simeq \left|
\begin{pmatrix} I&0\\0&V^*\end{pmatrix}\begin{pmatrix}
A&N\\0&B\end{pmatrix}\begin{pmatrix}
U^*&0\\0&I\end{pmatrix}\right|=\left|\begin{pmatrix}
|A^*|&N\\0&|B|\end{pmatrix}\right|
 $$
 and apply Theorem 2.1. \qquad $\Box$

 \vskip 10pt\noindent
 Using polar decompositions as above we may also obtain:

 \vskip 15pt\noindent
 {\bf Corollary 2.4.} {\it Let ${\mathbb{S}}$ be a  block-matrix of the
 form
 $$
 {\mathbb{S}}=\begin{pmatrix} A&I\\I&B\end{pmatrix}.
 $$
 in which $I$ is the identity. Let $f:[0,\infty)\to[0,\infty)$ be  concave and e-convex. Then, for all symmetric
 norms,
 $$
\Vert f(|{\mathbb{S}}|)\Vert \le \Vert f(|A^*|)+ 2f(I)+f(|B|) \Vert.
 $$
 }

 \vskip 10pt The
assumption in Theorem 2.1 requiring normality of each block is
rather special. The next Theorem meets the simple requirement that
the full matrix is Hermitian.

 \vskip 15pt\noindent
 {\bf Theorem 2.5.}  {\it Let ${\mathbb{A}}=[A_{i,\,j}]$ be a  Hermitian matrix  partitioned in blocks of same size and let
 $f$ be  a non-negative concave, e-convex function on $[0,\infty)$.
  Then, for all symmetric norms,
$$
\left\|\, f(|{\mathbb{A}}|) \,\right\| \le \left\|\,\sum
f(|A_{i,\,j}|) \,\right\|.
$$
}

\vskip 10pt\noindent {\bf Proof.} The proof of Theorem 2.1 actually
shows that for a general block-matrix  ${\mathbb{A}}=(A_{i,\,j})$
partitioned in blocks of same size, we have
$$
 \begin{pmatrix}f(|{\mathbb{A}}|)&0\\ 0&f(|{\mathbb{A}}|)\end{pmatrix}\prec_{w}
 \begin{pmatrix} \sum f(|A_{i,\,j}^*|)&0\\ 0& \sum f(|A_{i,\,j}|)\end{pmatrix}.
$$
for all non-negative concave, e-convex function $f$. Therefore
$$
\left\|\, f(|{\mathbb{A}}|) \,\right\| \le \max\left\{\,
\left\|\,\sum f(|A_{i,\,j}^*|) \,\right\|\,;\,\left\|\,\sum
f(|A_{i,\,j}|) \,\right\|\,\right\}
$$
Assuming ${\mathbb{A}}$ Hermitian, we have $A_{i,\,j}^*=A_{j,\,i}$
and Theorem 2.5 follows. \qquad $\Box$

\vskip 10pt\noindent Thus we have

 \vskip 10pt\noindent
 {\bf Corollary 2.6.}  {\it For all  Hermitian matrices partitioned in four blocks of same size,
 all  symmetric norms and $0\le p\le1$,

 $$
\left\| \,\left|
\begin{pmatrix}A&B^* \\B&C \end{pmatrix}\right|^p \,\right\|
  \le \left\|\, |A|^p+|B|^p+|B^*|^p+|C|^p\,\right\|.
 $$
}

\vskip 5pt\noindent It seems natural to propose:
 \vskip 10pt\noindent
 {\bf Conjecture 3.}  Theorem 2.5 holds for all symmetric norms.
\vskip 10pt\noindent
 as well as a possible stronger statement
 \vskip 10pt\noindent
 {\bf Conjecture 4.} Theorem 2.5 holds for all normal matrices partitioned in blocks of same size and all symmetric norms.

\vskip 10pt Now we state a result for a general full matrix. However
we have to confine to the trace norm, and the result can be stated
as a trace inequality improving (3).

 \vskip 10pt\noindent
 {\bf Theorem 2.7.}  {\it
Let ${\mathbb{A}}=[A_{i,\,j}]$ be a matrix partitioned in blocks of
arbitrary size and let  $f$ be a continuous function on $[0,\infty)$
with $f(0)\ge 0$. Then, if $f(\sqrt{t})$ is concave,
\begin{equation*}
{\mathrm{Tr\,}}f(|{\mathbb{A}}|)\, \le \, \sum {\mathrm{Tr\,}}
f(|A_{i,\,j}|) .
\end{equation*}
}

\vskip 10pt\noindent A special case of this theorem, with
$f(t)=t^p$, $2\ge p$, is a simple inequality for the Schatten
$p$-norms, $2\ge p$, of a square matrix $A=(a_{i\,j})$,
$$
\Vert A \Vert_p \le \left(\sum |a_{i,\,j}|^p\right)^{1/p}.
$$
For $p>2$, the reverse inequality holds. These inequalities go back
to \cite{GoMa} if not sooner.
 In contrast with Theorem 2.7, we cannot  extend Corollary 2.6
  for $p$ running over $[0,2]$. Let us give a
simple example for a partitioned positive matrix, $p=2$ and the
operator norm.

\vskip 10pt\noindent
 {\bf Example 2.8.} Let
 $$
 Z=\begin{pmatrix}A&B \\B&C \end{pmatrix}
 $$
 where
 $$
 A=\begin{pmatrix}2&0 \\0&1 \end{pmatrix}, \qquad
 C=\begin{pmatrix}1&0 \\0&2 \end{pmatrix}, \qquad
 B=\begin{pmatrix}0&\sqrt{2} \\\sqrt{2} &0 \end{pmatrix}.
 $$
 Then $\Vert Z^2\Vert_{\infty}=6+\sqrt{35}>9=\Vert
 A^2+2B^2+C^2\Vert_{\infty}$.

\vskip 15pt We turn to the proof of Theorem 2.7. We start with an
elementary well-known Lemma. For $A,\,B\ge 0$ the majorisation
relation $A\prec B$, or $B\succ A$, means $A\prec_w B$
 and ${\mathrm{Tr\,}}A={\mathrm{Tr\,}}B$. If $f(t)$ is a convex function, not
 necessarily increasing, $A\prec B\Rightarrow f(A)\prec_w f(B)$.
 This property is crucial for the proof of Theorem 2.7.

\vskip 10pt\noindent
 {\bf Lemma 3.}  {\it Let ${\mathbb{B}}=[B_{i,\,j}]$  be a positive block-matrix, with $1\le i,\,j\le m$. Then
 $$
  {\mathbb{B}}\succ \begin{pmatrix} B_{1,\,1}&\ &\ \\
  \ &\ddots &\ \\
  \ &\ &B_{m,\,m} \\
 \end{pmatrix}.
 $$  }

\vskip 10pt\noindent
 {\bf Proof. } Let
 $$
  {\mathbb{B}}= \begin{pmatrix} B_{1,\,1}&C\ \\
  C^*&D
 \end{pmatrix}.
 $$
 and note that
$$
   \begin{pmatrix} B_{1,\,1}&0\ \\
  0&D
 \end{pmatrix}
 = \frac{1}{2}\begin{pmatrix} B_{1,\,1}&C\ \\
  C^*&D
 \end{pmatrix}
 +
  \frac{1}{2}\begin{pmatrix} B_{1,\,1}&-C\ \\
  -C^*&D
  \end{pmatrix}.
 $$
 Hence,
 $$
  {\mathbb{B}}\succ \begin{pmatrix} B_{1,\,1}&0\ \\
  0&D
 \end{pmatrix}.
 $$
 Repeating this process with $D$, we see that the block-diagonal part of
 ${\mathbb{B}}$ lies in the convex hull of matrices unitarily
 congruent to ${\mathbb{B}}$. This implies the majorisation
 relation. \qquad $\Box$

 \vskip 15pt\noindent
 {\bf Proof of Theorem 2.7.} We have
 $|{\mathbb{A}}|^2={\mathbb{A}}^*{\mathbb{A}}$, hence Lemma 3 yields
 $$
 |{\mathbb{A}}|^2 \succ \begin{pmatrix} \sum_j |A_{j,1}|^2&\ &\ \\
  \ &\ddots &\ \\
  \ &\ & \sum_j |A_{j,m}|^2\\
 \end{pmatrix}.
 $$
Therefore, given any convex function $g(t)$ on $[0\infty)$,
$$
 g\left(|{\mathbb{A}}|^2\right) \succ_w \begin{pmatrix} g\left(\sum_j |A_{j,1}|^2\right)&\ &\ \\
  \ &\ddots &\ \\
  \ &\ & g\left(\sum_j |A_{j,m}|^2\right)\\
 \end{pmatrix}
 $$
so that
\begin{equation}
{\mathrm{Tr\,}} g\left(|{\mathbb{A}}|^2\right) \ge
{\mathrm{Tr\,}}g\left(\sum_j |A_{j,1}|^2\right) +\cdots+
{\mathrm{Tr\,}}g\left(\sum_j |A_{j,m}|^2\right)
\end{equation}
Now, if we assume furthermore that $g(0)\le 0$, we have
\begin{equation}
{\mathrm{Tr\,}} g(A+B) \ge {\mathrm{Tr\,}} g(A) +{\mathrm{Tr\,}}
g(B)
\end{equation}
for all $A,\, B\ge 0$, as a consequence of Theorem 1.1b (or of
Rotfeld's inequality). Combining (5) and (6) with
$g(t)=-f(\sqrt{t})$ ends the proof. \qquad $\Box$

\vskip 15pt


\vskip 20pt\noindent

 Jean-Christophe Bourin

 jcbourin@@univ-fcomte.fr

 Laboratoire de math\'ematiques

 Universit\'e de Franche-Comt\'e

 25030 Besan\c{c}on

\end{document}